\documentclass[draftclsnofoot,onecolumn]{IEEEtran}

\usepackage[utf8]{inputenc}
\usepackage[cmex10]{amsmath}
\usepackage{amsfonts}
\usepackage{amssymb}
\usepackage[english]{babel}

\newcommand{\eqn}[1]{(\ref{eqn:#1})}

\DeclareMathAlphabet{\mathpzc}{OT1}{pzc}{m}{it}

\begin{document}

\title{Moments and Absolute Moments of the Normal Distribution}

\author{\IEEEauthorblockN{Andreas Winkelbauer}\\%
\IEEEauthorblockA{Institute of Telecommunications, Vienna University of Technology\\%
Gusshausstrasse 25/389, 1040 Vienna, Austria\\email:~andreas.winkelbauer@nt.tuwien.ac.at}%
}

\maketitle

\begin{abstract}
We present formulas for the (raw and central) moments and absolute moments of the normal distribution. We note that these results are not new, yet many textbooks miss out on at least some of them. Hence, we believe that it is worthwhile to collect these formulas and their derivations in these notes.
\end{abstract}

\section{Introduction}
\label{sec:i}

Let $X \sim \mathcal{N}(\mu, \sigma^2)$ be a normal (Gaussian) random variable (RV) with mean $\mu = \mathrm{E}\lbrace X \rbrace$ and variance $\sigma^2 = \mathrm{E}\lbrace X^2 \rbrace - \mu^2$ (here, $\mathrm{E}\lbrace \cdot \rbrace$ denotes expectation). In what follows, we give formulas and derivations for 
$\mathrm{E}\big\lbrace X^\nu \big\rbrace$, 
$\mathrm{E}\big\lbrace ( X - \mu )^\nu \big\rbrace$, 
$\mathrm{E}\big\lbrace \lvert X \rvert^\nu \big\rbrace$, and 
$\mathrm{E}\big\lbrace \lvert X - \mu \rvert^\nu \big\rbrace$, 
i.e., for the (raw) moments, the central moments, the (raw) absolute moments, and the central absolute moments. We note that the formulas we present hold for real-valued $\nu > -1$.

The remainder of this text is structured as follows: Section \ref{sec:p} deals with preliminaries and introduces notation, particularly regarding some special functions. In Section \ref{sec:r} we present the results; the corresponding derivations are given in Section \ref{sec:d}.

\section{Preliminaries}
\label{sec:p}

We denote the standard deviation by $\sigma = \sqrt{\sigma^2}$. The imaginary unit is $j = \sqrt{-1}$ and $z^*$ denotes the complex conjugate of $z$. The nonnegative integers are denoted by $\mathbb{N}_0 = \mathbb{N} \cup \lbrace 0 \rbrace$. Next, we give the definitions of subsequently used special functions (cf.~\cite{specialFunc66}).
\begin{itemize}
\item \emph{Gamma function:}
\begin{equation}
\Gamma(z) \triangleq \int_0^\infty t^{z-1} e^{-t} dt.
\end{equation}
\item \emph{Rising factorial:}
\begin{align}
z^{\overline{n}} &\triangleq \frac{\Gamma(z+n)}{\Gamma(z)}\\
 &= z(z+1)\cdots(z+n-1), \quad n \in \mathbb{N}_0.
\end{align}
\item \emph{Double factorial:}
\begin{align}
z!! &\triangleq \sqrt{\frac{2^{z+1}}{\pi}}\, \Gamma\!\left(\frac{z}{2} + 1 \right)\\
 &= z \cdot (z-2) \cdot \ldots \cdot 3 \cdot 1, \quad z \in \mathbb{N} \text{ odd}.
\end{align}
\item \emph{Kummer's confluent hypergeometric functions:}
\begin{equation}
\Phi(\alpha, \gamma; z) \triangleq M(\alpha, \gamma, z) = {_1}F_1(\alpha; \gamma; z) = \sum_{n=0}^\infty \frac{\alpha^{\overline{n}}}{\gamma^{\overline{n}}} \frac{z^n}{n!}.
\end{equation}
\item \emph{Tricomi's confluent hypergeometric functions:}
\begin{equation}
\Psi(\alpha, \gamma; z) \triangleq U(\alpha, \gamma, z) = \frac{\Gamma(1-\gamma)}{\Gamma(\alpha-\gamma+1)} \Phi(\alpha, \gamma; z) + \frac{\Gamma(\gamma-1)}{\Gamma(\alpha)} z^{1-\gamma} \Phi(\alpha - \gamma + 1, 2 - \gamma; z).
\end{equation}
\item \emph{Parabolic cylinder functions:}
\begin{equation}
D_\nu(z) \triangleq 2^{\nu/2} e^{-z^2/4} \left[ \frac{\sqrt{\pi}}{\Gamma\!\left( \frac{1-\nu}{2} \right)} \Phi\!\left(-\frac{\nu}{2}, \frac{1}{2}; \frac{z^2}{2}\right) - \frac{\sqrt{2\pi}z}{\Gamma\!\left(-\frac{\nu}{2}\right)} \Phi\!\left(\frac{1-\nu}{2}, \frac{3}{2}; \frac{z^2}{2}\right) \right].
\end{equation}
\end{itemize}

\section{Results}
\label{sec:r}

In this section we give formulas for the raw/central (absolute) moments of a normal RV. If not noted otherwise, these results hold for $\nu > -1$.

\begin{itemize}
\item \emph{Raw moments:}
\begin{align}
\label{eqn:raw}
\mathrm{E}\big\lbrace X^\nu \big\rbrace & = (j\sigma)^\nu \exp\!\left(-\frac{\mu^2}{4\sigma^2}\right) D_\nu\!\left(-j\frac{\mu}{\sigma} \right)\\
 &= (j\sigma)^\nu 2^{\nu/2} \left[ \frac{\sqrt{\pi}}{\Gamma\!\left( \frac{1-\nu}{2} \right)} \Phi\!\left(-\frac{\nu}{2}, \frac{1}{2}; -\frac{\mu^2}{2\sigma^2}\right) + j \frac{\mu}{\sigma} \frac{\sqrt{2\pi}}{\Gamma\!\left( -\frac{\nu}{2} \right)} \Phi\!\left(\frac{1-\nu}{2}, \frac{3}{2}; -\frac{\mu^2}{2\sigma^2}\right) \right]\\
 &= (j\sigma)^\nu 2^{\nu/2} \cdot \left\lbrace \begin{array}{ll}
\Psi\!\left(-\frac{\nu}{2}, \frac{1}{2}; -\frac{\mu^2}{2\sigma^2}\right), & \mu \leq 0\\
\Psi^*\!\left(-\frac{\nu}{2}, \frac{1}{2}; -\frac{\mu^2}{2\sigma^2}\right), & \mu > 0
\end{array}\right.\\
 &= \left\lbrace \begin{array}{ll}
\sigma^\nu 2^{\nu/2} \frac{\Gamma\left( \frac{\nu+1}{2} \right)}{\sqrt{\pi}} \Phi\!\left(-\frac{\nu}{2}, \frac{1}{2}; -\frac{\mu^2}{2\sigma^2}\right), & \nu \in \mathbb{N}_0 \text{ even}\\
\mu \sigma^{\nu-1} 2^{(\nu+1)/2} \frac{\Gamma\left( \frac{\nu}{2} + 1 \right)}{\sqrt{\pi}} \Phi\!\left(\frac{1-\nu}{2}, \frac{3}{2}; -\frac{\mu^2}{2\sigma^2}\right), & \nu \in \mathbb{N}_0 \text{ odd}
\end{array}\right..
\end{align}
\item \emph{Central moments:}
\begin{align}
\label{eqn:central}
\mathrm{E}\big\lbrace ( X - \mu )^\nu \big\rbrace &= (j\sigma)^\nu 2^{\nu/2} \frac{\sqrt{\pi}}{\Gamma\left( \frac{1-\nu}{2} \right)}\\
\label{eqn:central2}
 &= (j\sigma)^\nu 2^{\nu/2} \cos( \pi \nu / 2 ) \frac{\Gamma\left( \frac{\nu+1}{2} \right)}{\sqrt{\pi}}\\
\label{eqn:central3}
 &= \big( 1 + (-1)^\nu \big) \sigma^\nu 2^{\nu/2 - 1} \frac{\Gamma\!\left( \frac{\nu+1}{2} \right)}{\sqrt{\pi}}\\
 &= \left\lbrace \begin{array}{ll}
\displaystyle \sigma^\nu (\nu-1)!!, & \nu \in \mathbb{N}_0 \text{ even}\\
\displaystyle 0, & \nu \in \mathbb{N}_0 \text{ odd}
\end{array}\right..
\end{align}
\item \emph{Raw absolute moments:}
\begin{align}
\label{eqn:rawabs}
\mathrm{E}\big\lbrace \lvert X \rvert^\nu \big\rbrace &= \sigma^\nu 2^{\nu/2} \frac{\Gamma\!\left( \frac{\nu+1}{2} \right)}{\sqrt{\pi}} \Phi\!\left( -\frac{\nu}{2}, \frac{1}{2}; - \frac{\mu^2}{2\sigma^2} \right).
\end{align}
\item \emph{Central absolute moments:}
\begin{align}
\label{eqn:centralabs}
\mathrm{E}\big\lbrace \lvert X - \mu \rvert^\nu \big\rbrace &= \sigma^\nu 2^{\nu/2} \frac{\Gamma\!\left( \frac{\nu+1}{2} \right)}{\sqrt{\pi}}.
\end{align}
\end{itemize}

\section{Derivations}
\label{sec:d}

In this section we give derivations for the previously presented results. Below we use the following two identities which hold for $\gamma \in \mathbb{R}, \nu > -1$ (cf.~\cite[Sec.~3.462]{gradshteyn6ed}):
\begin{align}
\label{eqn:hint1}
\int_{-\infty}^\infty (-jx)^\nu e^{-x^2+jx\gamma} dx &= \sqrt{2^{-\nu}\pi} e^{-\gamma^2/8} D_\nu\!\left( \frac{\gamma}{\sqrt{2}} \right),\\
\label{eqn:hint2}
\int_0^\infty x^\nu e^{-x^2-x\gamma} dx &= 2^{-(\nu+1)/2} \Gamma(\nu+1) e^{\gamma^2/8} D_{-\nu-1}\!\left( \frac{\gamma}{\sqrt{2}} \right).
\end{align}

\begin{itemize}
\item \emph{Raw moments \eqn{raw}:}
\begin{align}
\mathrm{E}\big\lbrace X^\nu \big\rbrace &= \frac{1}{\sqrt{2\pi\sigma^2}} \int_{-\infty}^{\infty} x^\nu \exp\!\left(-\frac{1}{2\sigma^2}(x-\mu)^2\right) dx\\
 &= \sqrt{\frac{2^\nu \sigma^{2\nu}}{\pi}} \exp\!\left(-\frac{\mu^2}{2\sigma^2}\right) \int_{-\infty}^{\infty} x^\nu \exp\!\left(-x^2 + x \frac{\mu}{\sigma} \sqrt{2} \right) dx\\
 &\hspace*{-1.1mm}\stackrel{\eqn{hint1}}{=} (j\sigma)^\nu \exp\!\left(-\frac{\mu^2}{4\sigma^2}\right) D_\nu\!\left(-j\frac{\mu}{\sigma} \right).
\end{align}
\item \emph{Central moments \eqn{central}:} Follows directly from \eqn{raw} with $\Phi\!\left(\alpha, \gamma; 0\right) = 1$ and, hence,
\begin{align}
D_\nu(0) = 2^{\nu/2} \frac{\sqrt{\pi}}{\Gamma\!\left(\frac{1-\nu}{2}\right)}.
\end{align}
To obtain \eqn{central2} from \eqn{central} we use the identity \cite[Sec.~8.334]{gradshteyn6ed}
\begin{align}
\Gamma\!\left(\frac{1+\nu}{2}\right) \Gamma\!\left(\frac{1-\nu}{2}\right) = \frac{\pi}{\cos( \pi \nu / 2 )}.
\end{align}
Then \eqn{central3} follows from \eqn{central2} by noting that
\begin{align}
\cos( \pi \nu / 2 ) = \frac{1 + \exp\!\left(j\pi\nu\right)}{2 \exp\!\left(j\pi\nu / 2\right)} = \frac{1 + (-1)^\nu}{2 j^\nu}.
\end{align}
\item \emph{Raw absolute moments \eqn{rawabs}:}
\begin{align}
\mathrm{E}\big\lbrace \lvert X \rvert^\nu \big\rbrace &= \frac{1}{\sqrt{2\pi\sigma^2}} \int_{-\infty}^{\infty} \lvert x \rvert^\nu \exp\!\left(-\frac{1}{2\sigma^2}(x-\mu)^2\right) dx\\
 &= \sqrt{\frac{2^\nu\sigma^{2\nu}}{\pi}}\, \exp\!\left(\!-\frac{\mu^2}{2\sigma^2} \right) \left[ \int_0^{\infty} x^\nu \exp\!\left(\!-x^2 -x \frac{\mu}{\sigma} \sqrt{2} \right) dx + \int_0^{\infty} x^\nu \exp\!\left(\!-x^2 +x \frac{\mu}{\sigma} \sqrt{2} \right) dx \right]\\
 &\hspace*{-1.1mm}\stackrel{\eqn{hint2}}{=} \sqrt{\frac{2^\nu\sigma^{2\nu}}{\pi}}\, \exp\!\left(\!-\frac{\mu^2}{4\sigma^2} \right) 2^{-(\nu+1)/2} \Gamma(\nu+1) \big( D_{-\nu-1}\!\left( \mu / \sigma \right) + D_{-\nu-1}\!\left(-\mu / \sigma \right) \big)\\
 &= \sqrt{\frac{\sigma^{2\nu}}{2^\nu}}\, \exp\!\left(\!-\frac{\mu^2}{2\sigma^2} \right) \frac{\Gamma(\nu+1)}{\Gamma(\nu/2+1)}\, \Phi\!\left( \frac{\nu+1}{2}, \frac{1}{2}; \frac{\mu^2}{2\sigma^2} \right)\\
 &= \sqrt{\frac{2^\nu\sigma^{2\nu}}{\pi}}\, \exp\!\left(\!-\frac{\mu^2}{2\sigma^2} \right) \Gamma\!\left( \frac{\nu+1}{2} \right) \Phi\!\left( \frac{\nu+1}{2}, \frac{1}{2}; \frac{\mu^2}{2\sigma^2} \right)\\
 &= \sigma^\nu 2^{\nu/2} \frac{\Gamma\!\left( \frac{\nu+1}{2} \right)}{\sqrt{\pi}} \Phi\!\left( -\frac{\nu}{2}, \frac{1}{2}; -\frac{\mu^2}{2\sigma^2} \right),
\end{align}
where we have used Kummer's transformation \cite[Sec.~9.212]{gradshteyn6ed}, i.e.,
\begin{align}
\Phi\!\left(\alpha, \gamma; z\right) = e^z \Phi\!\left(\gamma-\alpha, \gamma; -z\right),
\end{align}
in the last step.
\item \emph{Central absolute moments \eqn{centralabs}:} Follows directly from \eqn{rawabs} with $\Phi\!\left(\alpha, \gamma; 0\right) = 1$.
\end{itemize}

\renewcommand{\baselinestretch}{1.00}\normalsize\small
\bibliographystyle{IEEEtran}
\bibliography{tf-zentral,work_andreas}

\end{document}